\documentclass[12pt]{article}

\usepackage{latexsym,amsmath,amsthm,amssymb,fullpage}

\usepackage{color}

\usepackage{latexsym,amsmath,amsthm,amssymb,fullpage}

\usepackage[utf8]{inputenc}

\usepackage{amssymb}
\usepackage{epstopdf}
\usepackage{amsmath}
\usepackage{amsthm}
\usepackage{amsfonts}
\usepackage{enumerate}
\usepackage{epstopdf}

\newtheorem{theorem}{Theorem}
\newtheorem{lemma}[theorem]{Lemma}

\newtheorem{corollary}[theorem]{Corollary}

\newtheorem{proposition}[theorem]{Proposition}
\theoremstyle{definition}

\theoremstyle{remark}
\newtheorem{remark}[theorem]{Remark}

\newcommand\beq{\begin {equation}}
\newcommand\eeq{\end {equation}}
\newcommand\beqs{\begin {equation*}}
\newcommand\eeqs{\end {equation*}}

\newcommand\N{{\mathbb N}}

\newcommand\R{{\mathbb R}}

\newcommand\PP{{\mathbb P}}
\newcommand\E{{\mathbb E}}
\newcommand\sph{{\mathbb S}}

\begin {document}

$$   $$

\vskip 4mm
\begin{center}
{\LARGE \bf Some results on the optimal matching}
\vskip 2mm
{\LARGE \bf problem for the Jacobi model}
\vskip 10mm
\textit{\large Jie-Xiang Zhu} \footnote{Research supported by China Scholarship Council (No. 201806100097) and partially supported by NSF of China (Grants No. 11625102 and No. 11571077)}\\
\vskip 3mm
{Fudan University, People's Republic of China, \\and University of Toulouse, France}\\
\end{center}
\vskip 7mm

\vskip 4mm

\begin {abstract}
We establish some exact asymptotic results for a matching problem
with respect to a family of beta distributions. Let $X_1, \ldots, X_n$ be independent random variables with common distribution the symmetric Jacobi measure
$d\mu (x) = C_d (1-x^2)^{\frac d2 -1} dx$ with dimension $ d \geq 1$ on $[-1, 1]$,
and let $\mu_n = \frac{1}{n} \sum_{i = 1}^{n} \delta_{X_i}$ be the associated empirical measure.
We show that
\begin{align*}
\lim_{n \to \infty} n\E \left[ W_2^2( \mu^n, \mu ) \right] = \sum_{k = 1}^{\infty} \frac{1}{k(k+d-1)},
\end{align*}
where $W_2$ is the quadratic Kantorovich distance with respect to the intrinsic cost
$\rho(x, y) = |\arccos(x) - \arccos (y)|$, $(x, y) \in [-1, 1]^2$, associated to the model.
When $\mu$ is the product of two Jacobi measures with dimensions $d$ and
$d'$ respectively, then
$$
\E \left[ W_2^2( \mu^n, \mu ) \right] \approx  \frac{\log n}{n} \, .
$$
In the particular case $d = d' = 1$ (corresponding to the product of arcsine laws),
\begin{align*}
\lim_{n \to \infty} \frac{n}{\log n} \E \left[ W_2^2( \mu^n, \mu ) \right] = \frac{\pi}{4}.
\end{align*}
Similar results do hold for non-symmetric Jacobi distributions.
The proofs are based on the recent PDE and mass transportation approach developed
by L.~Ambrosio, F.~Stra and D.~Trevisan.
\end {abstract}



\section{ Introduction and main results} \label{S1}
\setcounter{equation}{0}
Optimal matching problems are very classical in computer science, physics and mathematics,
and have been widely investigated from different viewpoints. One general formulation, of particular
interest in probability theory and mathematical statistics, is expressed in terms
of Kantorovich (Wasserstein)
distances. Let $(M, \rho)$ be a
(complete separable) metric space. Given $p \geq 1$, the Kantorovich distance (see e.g.~\cite{V09}) between two probability measures $\nu$ and $\mu$
on the Borel sets of $M$ with a finite $p$-th moment is defined by
\begin{align*}
W_p(\nu, \mu) = \inf \left( \int_{M \times M} \rho^p (x, y) d\pi(x, y) \right)^{\frac{1}{p}}
\end{align*}
where the infimum is taken over all couplings $\pi$ on $M \times M$ with respective marginals $\nu$
and $\mu$. In this work, we mainly focus on the quadratic one ($p = 2$).

Given then a probability measure $\mu$ on the Borel sets of $(M, \rho)$, denote by $X_1, \ldots, X_n, \ldots$ independent random variables with common distribution $\mu$, and let
$\mu^n = \frac{1}{n}\sum_{i=1}^{n} \delta_{X_i}$, $n \geq 1$, be the empirical measure
on the sample $(X_1, \ldots, X_n)$. One version of the optimal matching problem, at a first order,
is to estimate the order of growth, and possibly the renormalized limit, of
$$
\E \left[ W_2^2(\mu^n, \mu) \right]
$$
as $n \to \infty$.

A reasonable answer to this question involves conditions on the nature of the metric space $(M,\rho)$
and on the distribution $\mu$. In the past three decades, several authors made contributions to this problem,
and we mention here a few relevant ones. One first well-known result is the
Ajtai-Koml\'{o}s-Tusn\'{a}dy theorem \cite{AKT84}
expressing that for $\mu$ the uniform distribution on the unit cube $[0, 1]^k$ in $\R^k$, $k \geq 1$,
\begin{align} \label {AKT}
\E \left[ W_p^p( \mu^n, \mu ) \right] \approx \left\{ \begin{array}{cc}
\frac{1}{n^{p/2}},  & k = 1, \\
\Big ( \frac{\log n}{n} \Big)^{p/2},  & k = 2, \\
\frac{1}{n^{p/k}}, & k \geq 3.
\end{array}
\right.
\end{align}

In this paper we use the notation $A \lesssim B$ if there exists a constant $C > 0$
(depending on the underlying model but not on $n$)
such that $A \leq C B$. We also write $A \approx B$, if $A \lesssim B$ and $B \lesssim A$.
(It will also be implicit throughout the paper that $n \geq 1$ is large enough so that all
the relevant quantities involving $n$ are well defined.)
This first result was completed, refined and extended to more general situations and examples
in several works. We refer to the monograph \cite {T14} by M.~Talagrand
for an account and relevant references on this result and its
developments (see also \cite {Y98,HS92}).
Additional recent advances include \cite {BM02,BB13,DSS13,FG15,HPZ18} to cite a few.
For an extensive discussion of the one-dimensional optimal matching problem, see \cite {BL16}.

It should be emphasized that in \eqref {AKT} the most delicate case is $k=2$ since the rate
does not reflect the uniform spacings in $[0,1]^2$ and
both local and global discrepancies have to be combined (cf.~\cite {T14}). This is also
the case for $k=1$, but here it is classical that the optimal matching is achieved by
monotone rearrangement. In particular, order statistics may be used to produce the exact expression
\begin{align} \label {uniform}
\E \left[ W_2^2( \mu^n, \mu ) \right] = \frac{1}{6n}
\end{align}
for $\mu$ the uniform distribution on $[0, 1]$ (cf. e.g.~\cite {BL16}).
Order statistics can also be applied to estimate the mean Kantorovich distance
for log-concave distributions on the real line (see Section~6 of \cite{BL16}). For
the purpose of this work, let us observe in particular that if
$d\mu(x) = C_{d} (1-x^2)^{\frac{d}{2}-1} dx$, $d \geq 1$, is the symmetric Jacobi
(beta) distribution on $[-1,1]$
(for a detailed definition, see Section~\ref{S2}), then
\begin{align} \label {jacobieuclidean}
\E \left[ W_2^2( \mu^n, \mu ) \right] \approx \frac{1}{n}.
\end{align}

Recently, a major achievement was achieved by L.~Ambrosio, F.~Stra and D.~Trevisan \cite{AST19}
who answered rigorously a conjecture put forward in \cite {CLPS14}. They indeed provided the exact value
of the limiting behaviour of $\E [ W_2^2( \mu^n, \mu ) ]$
in the Ajtai-Koml\'{o}s-Tusn\'{a}dy theorem for $k=2$, namely
\begin{align} \label {AST}
\lim_{n \to \infty} \frac{n}{\log n} \E\left[ W_2^2( \mu^n, \mu ) \right] = \frac{1}{4\pi}.
\end{align}
This is actually one of the very few explicit limits known in this setting,
even in dimension one actually
(some unknown limits are achieved via subadditivity in higher dimensions in \cite {DY95,BM02,BB13,DSS13}).

More generally, based on the PDE ansatz of \cite {CLPS14}, the authors of
\cite {AST19} proved that the same limit
\eqref {AST} holds true for $\mu$ the normalized uniform distribution
on a compact 2-dimensional Riemannian manifold $(M,g)$ without boundary.
The validity of the PDE ansatz of \cite {CLPS14}
is based on a fine analysis of a linearized Monge-Amp\`{e}re equation and mass transportation tools, which highly depend on properties of the associated heat kernel $p_t$. It is expected that the method can be further developed for the case $d \geq 3$, but this is mostly conjectural at this point.

\bigskip

Relying on the methodology emphasized in \cite {AST19},
the purpose of this work is to provide some further examples, in (topological) dimensions $1$ and $2$,
where the exact asymptotic behaviour of $\E [ W_2^2( \mu^n, \mu )]$ may be achieved.
These examples are the aforementioned symmetric Jacobi measures
$$
d\mu(x) = C_{d} (1-x^2)^{\frac{d}{2}-1} dx \quad \mbox {on} \, \, [-1,1],
$$
(also called beta laws). Note that when $d=2$, this is the uniform
measure on $[-1,1]$ while when $d=1$, $\mu$ is the famous arcsine law.
However, to fully develop the PDE
and optimal transport approach, the interval $[-1,1]$, and thus the Kantorovich distance
$W_2$, will be equipped with the intrinsic
distance $\rho(x, y) = |\arccos(x) - \arccos (y)|$, $(x, y) \in [-1,1]^2$, of the underlying Jacobi
model (see Section~\ref{S2} below). As indeed described below, the intrinsic distance is inherited
from the geodesic metric on the sphere (Jacobi measures are projections of the spherical uniform measure on a diameter when $d$ is an integer),
and compares to the Euclidean one, but for the limit to hold via the PDE approach,
the intrinsic distance is in force. In particular, Theorem~\ref {t11} below for $d=2$
(for which the limit is $1$) is not equivalent to \eqref {uniform}.

In this setting, the main results of this work read as follows.

\begin{theorem} \label{t11}
 If $\mu$ is the symmetric Jacobi measure with $d \geq 1$, then
\begin{align*}
\lim_{n \to \infty} n\E \left[ W_2^2( \mu^n, \mu ) \right] = \sum_{k = 1}^{\infty} \frac{1}{k(k+d-1)}.
\end{align*}
\end{theorem}

Note that when $d \geq 2$ is an integer,
$$
\sum_{k = 1}^{\infty} \frac{1}{k(k+d-1)} = \frac {1}{d-1} \sum_{k=1}^{d-1} \frac 1k
$$
which is of the order $\frac {\log d}{d}$ as $d \to \infty$.

\begin{theorem} \label{t12}
If $\mu$ is the product measure of two symmetric Jacobi measures with dimensions
$d, d' \geq 1$ respectively, then
$$
\E \left[ W_2^2( \mu^n, \mu ) \right] \approx \frac{\log n}{n}
$$
(where $W_2$ is defined via the product distance of the intrinsic metric $\rho$).
\end{theorem}

\begin{theorem} \label{t13}
If $\mu$ is the product measure of two symmetric Jacobi measures with dimensions $d = d' = 1$
(i.e. arcsine laws), then
\begin{align*}
\lim_{n \to \infty} \frac{n}{\log n} \, \E \left[ W_2^2( \mu^n, \mu ) \right] = \frac{\pi}{4}.
\end{align*}
\end{theorem}

It is likely that Theorem~\ref {t13} holds true when $d+ d' =2$, but the probabilistic
tail estimates require $d+d' \leq 2$ while at the same time $d, d'\geq 1$ is requested
for the heat kernel bounds.
If $\mu$ is the product measure of $k \geq 3$ symmetric Jacobi measures with dimensions $d_j \geq 1$,
$j = 1, \ldots, k $ respectively, Theorem~\ref {t12} extends in the form
\beq \label {eq.kjacobi}
\E \left[ W_2^2( \mu^n, \mu ) \right] \approx  \frac{1}{n^{2/k}}
\eeq
in accordance with \eqref {AKT}.

The results may also be formulated in the bipartite form, corresponding really to
the optimal matching problem. If $X_1, \ldots, X_n$ and $Y_1, \ldots, Y_n$ independent random variables with common distribution $\mu$, let $\mu^n = \frac{1}{n}\sum_{i=1}^{n} \delta_{X_i}$
and  $\nu^n = \frac{1}{n}\sum_{i=1}^{n} \delta_{Y_i}$ be the corresponding
empirical measures. Then, if $\mu$ is the symmetric Jacobi measure with $d \geq 1$,
\beq \label {bipartite1}
\lim_{n \to \infty} n\E \left[ W_2^2( \mu^n, \nu^n ) \right] = \sum_{k = 1}^{\infty} \frac{2}{k(k+d-1)}
\eeq
and if $\mu$ is the product measure of two symmetric Jacobi measures with
dimensions $d = d' = 1$,
\beq \label {bipartite2}
\lim_{n \to \infty} \frac{n}{\log n} \E \left[ W_2^2( \mu^n, \nu^n ) \right] = \frac{\pi}{2}.
\eeq

As discussed during the proofs below, the main conclusions above actually also apply
to non-symmetric Jacobi distributions
$$
d \mu (x) = C_{\alpha, \beta} (1 - x)^{\alpha - 1}(1 + x)^{\beta - 1}dx
$$
with $\alpha, \beta \geq \frac{1}{2}$. In this case, Theorem~\ref {t11} takes the same form with
$d = \alpha + \beta$. In Theorem~\ref {t12}, $\mu$ may be taken
as the product of two such non-symmetric Jacobi measures with $\alpha, \beta \geq \frac{1}{2}$,
$\alpha', \beta' \geq \frac{1}{2}$. Theorem~\ref {t13} is unchanged
since the restrictions $\alpha, \beta \geq \frac{1}{2}$ imply that we are
necessarily in the symmetric case.

As mentioned above, the symmetric Jacobi model when $d$ is an integer
can be seen as a quotient of the $d$-dimensional sphere $\sph^d$,
and as a Markov Triple (see \cite {BGL14} and Section~\ref {S2} below)
satisfies a curvature-dimension condition $CD(d-1,d)$ as the sphere.
The scaling and limit in Theorem~\ref {t11} however do not reflect this property. This is due
to the fact that the PDE approach is a tradeoff between the (small time) heat kernel behaviour
$$
p_s(x, x) \leq Cs^{- \frac{d}{2}}, \quad 0 < s \leq 1,
$$
which is alike $d$-dimensional for the Jacobi semigroup, and the trace estimate
\begin{align} \label{tr}
\int_{-1}^{1} p_s(x, x) d\mu(x) \leq Cs^{-\frac{1}{2}}, \quad 0 < s \leq 1,
\end{align}
which actually reflects a one-dimensional feature. While the heat kernel and trace
small time asymptotics are of the same order
on a compact Riemannian manifold, the Jacobi model is intermediate in this respect.

This paper is therefore a contribution towards further examples where
exact limiting behaviours in the optimal matching problem can be achieved, following the conjectures
in \cite {CLPS14}. At this point very few results are available on this rather delicate problem.
The methodology developed to this task is the PDE and optimal transport approach of \cite {AST19},
that is followed quite closely. In particular, some intermediate results and arguments of \cite {AST19}
may be immediately adapted to the Jacobi setting, so that we only quote them without detailed proofs.
In addition, the investigation also benefits from the recent contribution \cite {AG18}, in particular at the level
of a strengthened regularization step. At the same time, the full organization of the proof
presented here produces several simplified arguments and steps, which may be used in return
to streamline the proofs in \cite {AST19,AG18}. In particular, Riesz transform bounds are not used.

Turning to the content of the paper,
Section~\ref {S2} presents the fundamental Markov Triple structure from \cite {BGL14},
the example of the Jacobi model, and several properties of the associated (Jacobi) semigroup. From
Section~\ref {S3} to Section~\ref {S5}, the main technical tools are introduced, including density fluctuation bounds, energy estimates and a refined contractivity estimate, following \cite {AST19} and \cite {AG18}.
In Section~\ref {S6}, we prove Theorem~\ref{t11} by treating upper and lower bounds respectively, according to the same outline as in \cite{AST19}. The proof of Theorem~\ref {t13} is similar. In Section~\ref {S7},
we prove Theorem~\ref{t12} by means of Proposition~3 of \cite{L17} and a simple comparison argument.

\section{Markov Triple, Jacobi model and properties of the semigroup} \label{S2}
\setcounter{equation}{0}

This section is a brief exposition of the Jacobi model within the setting of Markov Triples
as put forward in \cite {BGL14} (see in particular therein Chapters 1 to 3).
To start with, we briefly recall the basic definition of a Markov Triple $(E, \mu, \Gamma)$
referring to \cite {BGL14} for the complete picture.

On a measurable space $(E, \mathcal{F})$, a Markov semigroup $(P_t)_{t \geq 0}$ is a family of operators defined on some subset of real-valued measurable functions on $(E, \mathcal{F})$ satisfying the positivity and mass conservation properties. A $\sigma$-finite measure $\mu$ is said to be invariant
for the semigroup $(P_t)_{t \geq 0}$ if for every bounded positive measurable function $f$ on $E$ and $t \geq 0$,
\begin{align*}
\int_{E} P_t f d\mu = \int_{E} f d\mu.
\end{align*}

Denote the infinitesimal generator of $(P_t)_{t \geq 0}$ by $L$. The associated carr\'{e} du champ operator $\Gamma$ is defined by
\begin{align*}
\Gamma(f, g) = \frac{1}{2}\left[ L(fg) - fL(g) - gL(f) \right],
\end{align*}
where $f$ and $g$ are elements of a suitable algebra $\mathcal {A}$ of
(smooth) functions. The integration by parts formula expresses in this context that
\begin{align*}
\int_{E} g L(f) d\mu = - \int_{E} \Gamma(f, g) d\mu.
\end{align*}
In sequel, $\Gamma(f, f)$ is abbreviated as $\Gamma(f)$.
The family $(E, \mu, \Gamma)$ is then called a Markov Triple.

The Markov Triple structure induces in addition a notion of curvature-dimension condition
$CD(K,N)$, $K \in \R$, $N \geq 1$, which extends the Ricci curvature lower bounds
in Riemannian manifolds.

As an example, let $E$ be $\R^k$, and let $\mu$ be Lebesgue measure. For the
carr\'e du champ operator $\Gamma(f, g) = \sum_{i=1}^{k} \partial_i f \partial_i g$
on smooth functions $f,g : \R^k \to \R$, $(\R^k, \mu, \Gamma)$ is a Markov Triple, the
generator $L$ being the standard Laplace operator $\Delta $ on $\R^k$. More generally,
a (smooth complete) Riemannian manifold with the Riemannian volume element and the
Riemannian length squared $\Gamma (f) = |\nabla f|^2$ of the gradient of a smooth function
$f$ defines similarly a Markov Triple. An $k$-dimensional manifold $(M,g)$ with Ricci curvature
bounded below by $K$ satisfies the curvature-dimension condition $CD(K,k)$.

The Markov Triple definition $(E, \mu, \Gamma)$ also
induces an intrinsic metric structure on $E$ by means of the distance
$$
\rho (x, y) = \mathrm {esssup} [f(x) - f(y)], \quad (x,y) \in E \times E,
$$
where the essential supremum is taken over all $1$-Lipschitz functions $f$
on $E$ such that ${\Gamma(f) \leq 1}$ $\mu$-almost everywhere. For example,
on the preceding example $(\R^k, \mu, \Gamma)$, $\rho$ is simply the Euclidean metric on $\R^k$.
When dealing with the Kantorovich distance $W_2$ on a Markov Triple $(E, \mu, \Gamma)$,
the underlying distance will be the intrinsic distance $\rho$ as just defined.

Markov Triples are extensively discussed and studied in the monograph \cite {BGL14}.
The Jacobi model is one family of examples.
On the interval $E = [-1, 1]$, the Jacobi measure with parameters
$\alpha , \beta >0$ is defined by
$$
d \mu_{\alpha, \beta}(x) = C_{\alpha, \beta} (1 - x)^{\alpha - 1}(1 + x)^{\beta - 1}dx
$$
where $C_{\alpha,\beta} > 0$ is a normalization constant and $dx$ is Lebesgue measure.
(While for $0 < \alpha < 1$ or $0 < \beta < 1$, $\mu_{\alpha, \beta}$ is formally only defined
on $(-1,1)$, the relevant Markov Triple properties emphasized below do hold similarly on the closed
interval $E = [-1, 1]$ that we thus keep for simplicity in the exposition.)
Throughout the investigation, it will be assumed that $\alpha, \beta \geq \frac 12$ to ensure
the necessary analytic properties towards the main conclusions -- specifically the ultracontractivity
property (UC), see below.

The associated Jacobi operator is acting on smooth functions $f$ as
$$
L_{\alpha, \beta}f = (1 - x^2)f''- [(\alpha + \beta)x + \alpha - \beta]f',
$$
and the carr\'{e} du champ operator is given by $\Gamma( f ) = (1 - x^2){f'}^2$.
The family $([-1, 1], \mu_{\alpha, \beta}, \Gamma)$ is therefore
a Markov Triple, called the Jacobi model (of parameters $\alpha, \beta$). It is easy to check that the intrinsic
distance $\rho$ between $x$ and $y \in [-1, 1]$ is given by
\begin{align*}
\rho (x, y) = \left| \int_{x}^{y} \frac{1}{\sqrt{1-u^2}} du \right| = | \arccos(x) - \arccos(y) |.
\end{align*}

When $\alpha = \beta = \frac d2$ for some $d \geq 1$,
we speak of the symmetric Jacobi model of dimension $d$, with thus the symmetric
Jacobi distribution $d\mu_{\frac{d}{2}, \frac{d}{2}}(x) = C_{d} (1-x^2)^{\frac{d}{2}-1} dx$.
The ``dimension" terminology has a geometric flavour. Indeed, when
$d \geq 1$ is an integer,
the symmetric Jacobi model can be seen as a quotient of the sphere $\sph^d$ in $\R^{d+1}$ via the projection
operator $T : x =(x_1, \ldots, x_{d+1}) \in \sph^d \mapsto x_1 \in [-1, 1] $.
If $X$ is distributed uniformly on $\sph^d$, then the law of $T(X)\in [-1 , 1]$ is precisely
$\mu_{\frac{d}{2}, \frac{d}{2}}$. Moreover, the Jacobi operator $L_{\frac{d}{2}, \frac{d}{2}}$
can be seen as an image of the Laplace-Beltrami operator
$\Delta_{\sph^d}$ on $\sph^d$ via its action on functions depending on only
one coordinate in $\R^{d+1}$. According to this geometric description,
the Markov Triple $([-1, 1], \mu_{\frac{d}{2}, \frac{d}{2}}, \Gamma)$
satisfies the curvature-dimension condition $CD(d-1,d)$ as the sphere $\sph^d$,
a property which extends to non-integer $d$.

In Theorems~\ref {t12} and \ref {t13},
we also consider the product Markov Triple generated by two or more Jacobi models.
On the square $E = [-1, 1]^2$, the product measure $\mu$ of two Jacobi
measures with parameters $\alpha, \beta >0$, $\alpha ', \beta' >0$ respectively, is defined by
$$
d\mu(x, y) = C (1-x)^{\alpha-1} (1+x)^{\beta -1} (1-y)^{\alpha'-1} (1+y)^{\beta' -1} dx dy
$$
where $C > 0$ is a normalization constant and $dx dy$ is Lebesgue measure.
The associated carr\'{e} du champ operator is
\begin{align*}
\Gamma( f )(x, y) =  (1 - x^2)(\partial_x f)^2 + (1 - y^2)(\partial_y f )^2.
\end{align*}
Then the intrinsic distance $\rho$ between $x = (x_1, x_2)$ and $y = (y_1, y_2) \in [-1, 1]^2$ satisfies
\begin{align*}
\rho(x , y) = \sqrt{ \rho_1^2(x_1, y_1) + \rho_2^2(x_2, y_2) },
\end{align*}
where $\rho_i$, $i = 1, 2$, are the intrinsic distances of the respective Jacobi Triples.
When $\alpha = \beta = \frac d2$
and $\alpha' = \beta' = \frac {d'}{2}$, the product Jacobi model is of
curvature-dimension $CD(\min (d,d') -1, d+d')$. Similar properties hold for the Markov triple generated by more than two Jacobi models.

To streamline the exposition and emphasize the main ideas, we concentrate on
the symmetric Jacobi model throughout the proofs and the arguments. The necessary
modifications to cover the non-symmetric case will be mentioned at the appropriate places.
Accordingly, in the sequel, we drop for simplicity
the subscripts $\frac{d}{2}$ in $\mu_{\frac{d}{2}, \frac{d}{2}}$, and omit the
word ``symmetric".

Given then a Jacobi Triple $([-1, 1], \mu, \Gamma)$
(where $ \mu = \mu_{\frac{d}{2}, \frac{d}{2}}$ with $d \geq 1$ fixed),
denote by $(P_t)_{t \geq 0}$ the symmetric Markov semigroup with infinitesimal generator $L$.
In following sections, we will investigate properties of the solution of the (Poisson) equation $Lf = u$.
To this task, note that, formally, $(-L)^{-1}$ is described as
\begin{align*}
(-L)^{-1} = \int_{0}^{\infty} P_t dt,
\end{align*}
acting on mean zero functions in the suitable domain. On the other hand, it is known that
for the Jacobi model, the eigenvalues of $-L$ are given by the sequence
$\lambda_k = k(k+d-1)$, $k \in \N$. The corresponding eigenvectors $J_k$
are the Jacobi polynomials of degree $k$ which form an orthogonal basis
of $L^2(\mu)$. In this spectral description, $(-L)^{-1}$ can then be represented by
\begin{align*}
(-L)^{-1} f = \sum_{k=1}^{\infty} \frac{1}{\lambda_k} f_k J_k,
\end{align*}
where $f_k = \int_{-1}^{1} f J_k d\mu$.

The operators $P_t$ of the Markov semigroup are represented by kernels
$$
P_t f(x) = \int_{-1}^1 f(y) p_t(x,y) d\mu (y), \quad t >0, \, \, x \in [-1,1],
$$
with respect to the invariant measure $\mu$, called the heat kernel. In the preceding
spectral representation,
$$
P_t f = \sum_{k=0}^{\infty} e^{-\lambda_kt} f_k J_k,
$$
and the trace formula is expressed by
$$
\int_{-1}^1 p_t (x,x) d\mu (x) =  \sum_{k=0}^{\infty} e^{-t\lambda_k} .
$$

The PDE and optimal transport
approach developed in \cite {AST19} is based on several quantitative properties of
the Markov semigroup $(P_t)_{t \geq 0}$ and the associated heat kernel $p_t(x,y)$, which we gather
in the following list.

There exist constants $C > 0$ and $K \geq 0$, possibly changing from line to line, such that:

\bigskip

(SG) Spectral gap: $\| P_t f \|_2 \leq e^{-Ct} \|  f \|_2$ for any $f$ with $\int_{-1}^{1} f d\mu = 0$
and $t\geq 0$.

(UC) Ultracontractivity: $p_t(x, y)  \leq C t^{-d/2}$ for any $x, y\in [-1, 1]$
and $0< t \leq 1$.

(GE) Gradient estimate: $\mathrm {Lip}(p_t(x, \cdot)) \leq C t^{-(d+1)/2}$ for any $x \in [-1, 1]$
and $0 < t \leq 1$.

(DR) Dispersion rate: $\int_{-1}^{1} \int_{-1}^{1} \rho^2(x, y) p_t(x, y) d\mu(x) d\mu (y) \leq Ct$,
$t > 0$.

(SGB) Strong gradient bound: $\sqrt{\Gamma(P_t f)} \leq e^{-Kt} P_t(\sqrt{\Gamma(f)})$
for any smooth $f$ and $t \geq 0$.

\bigskip

These properties are presented and detailed in \cite {AST19} in the context of a compact Riemannian
manifold (without boundary). We provide here the necessary arguments in case of the Jacobi model
under investigation, referring to \cite {BGL14} for the relevant properties needed to this task.

The spectral gap condition (SG) may be seen as a consequence of the spectral decomposition. Namely,
mean zero functions are orthogonal to constants, so that the exponential decay (SG) holds true
with $C = \lambda_1 > 0$. The other four properties actually derive from the curvature-dimension
condition $CD(d-1, d)$ of the Jacobi model.
Under this condition, a Sobolev inequality (of dimension $d$) holds, yielding equivalently the
uniform heat kernel bound (UC). (The case $d=1$ is a bit particular here since only $CD(0,1)$
is then available while the proof of the Sobolev inequality under a curvature-dimension condition
requires positive curvature. We may nevertheless rely then on the direct, and more precise,
heat kernel bounds developed in \cite {NS13}. This reference is also used to justify the restriction
$d \geq 1$ in this investigation since the heat kernel bounds do not seem to have been clearly put forward
in the range $0 < d < 1$.) The strong gradient bound (SGB), actually equivalent to
$CD(K, \infty)$, therefore holds for the Jacobi model with $K = d-1 \geq 0$.

Some more care has to be taken with (GE) and (DR). First, a curvature condition
$CD(0, \infty)$ implies the local Poincar\'{e} inequality
\begin{align*}
\Gamma(P_t f) \leq \frac{1}{2t} \big [ P_t(f^2) - (P_t f)^2 \big]
\end{align*}
for any (smooth) $f$ and $t >0$. For $s>0$ and $x \in [-1,1]$ fixed, the choice of
$f = p_s (x, \cdot)$ together with the semigroup property imply that
\begin{align*}
\| \Gamma( p_{t+s} (x, \cdot) ) \|_\infty
  \leq \frac{1}{2t} \, \| P_t(f^2) \|_\infty
   \leq \frac{1}{2t} \| P_t \|_{1 \to \infty }\int_{-1}^{1} f^2 d\mu
   = \frac{1}{2t} \| P_t \|_{1 \to \infty}p_{2s}(x, x).
\end{align*}
Now (UC) tells us that $\| P_t \|_{1 \to \infty}$ and $p_{2s}(x, x)$ are bounded from above by
$\frac{C}{t^{d/2}}$ and $\frac{C}{{2s}^{d/2}}$ respectively, so that (GE) follows by choosing $s = t$.

Next we turn to (DR). To this task we investigate, for any 1-Lipschitz function $f$, the integral
$$
\int_{-1}^{1} \int_{-1}^{1} (f(x) - f(y))^2 p_t(x, y) d\mu(x) d\mu (y).
$$
By the semigroup properties,
\begin{align*}
&\int_{-1}^{1} \int_{-1}^{1} (f(x) - f(y))^2 p_t(x, y) d\mu(x) d\mu (y)\\
& = 2 \left( \int_{-1}^{1}\int_{-1}^{1} f^2(x) p_t(x, y) d\mu(x) d\mu (y) - \int_{-1}^{1}\int_{-1}^{1} f(x) f(y) p_t(x, y) d\mu(x) d\mu (y) \right)\\
& = 2 \left( \int_{-1}^{1}P_t(f^2)(y)  d\mu (y) - \int_{-1}^{1} f(y) P_t(f)(y) d\mu(y) \right)\\
& = 2 \left( \int_{-1}^{1}f^2(y) d\mu(x) d\mu (y) - \int_{-1}^{1}P_{\frac{t}{2}}(f)^2(y)  d\mu(y) \right).
\end{align*}
Then we make use of $\Gamma$-calculus to obtain that
\begin{align*}
& \int_{-1}^{1}f^2(y) d\mu(x) d\mu (y) - \int_{-1}^{1}P_{\frac{t}{2}}(f)^2(y)  d\mu(y) \\
& = - \int_{0}^{\frac{t}{2}} \frac{d}{ds} \left( \int_{-1}^{1} P_s (f)^2(y) d\mu(y) \right) ds \\
& = -2 \int_{0}^{\frac{t}{2}} \left( \int_{-1}^{1} P_s (f)(y) LP_s(f)(y) d\mu(y) \right) ds \\
& = 2 \int_{0}^{\frac{t}{2}} \left( \int_{-1}^{1} \Gamma(P_s(f))(y) d\mu(y) \right) ds \\
& \leq 2 \int_{0}^{\frac{t}{2}} \left( \int_{-1}^{1} e^{-2Ks} P_s(\Gamma(f))(y) d\mu(y) \right) ds \leq 2 \int_{0}^{\frac{t}{2}}  e^{-2Ks}  ds.
\end{align*}
The last line comes from (SGB) and the fact that $\Gamma(f) \leq 1$. Hence for all $t > 0$,
\begin{align} \label{q1}
 \int_{-1}^{1} \int_{-1}^{1} (f(x) - f(y))^2 p_t(x, y) d\mu(x) d\mu (y)\leq \frac{2}{K}(1-e^{-Kt}) \leq 2t .
\end{align}

The preceding actually holds in a general context. But for the specific Jacobi model, we may choose
$f(x) = \arccos(x)$ in \eqref{q1} for which $\Gamma(\arccos(\cdot)) = 1$. Since
${\rho (x, y) = | \arccos(x) - \arccos(y) |}$, this is exactly (DR).

\begin {remark}
Using the sharp estimates on the Jacobi kernel $p_t$ from \cite{NS13},
it is possible to reach a stronger version of (DR) in the sense that for some $C>0$,
for all $x \in [-1, 1]$ and $t >0$, $\int_{-1}^{1} \rho^2(x, y) p_t(x, y) d\mu (y) \leq Ct$.
\end {remark}

As observed in \cite {AST19}, we record for further purposes that
(SGB) and (SG) also imply (with the same proof as in \cite {AST19}) that there exists a positive constant $C$
such that for every $f \in L^2(\mu)$ with $Lf \in L^{\infty} (\mu)$,
\begin{align} \label{i1}
\| \sqrt{\Gamma(f)} \|_{\infty} \leq C \| Lf \|_{\infty}.
\end{align}
Moreover, the preceding heat kernel properties are stable under product, and therefore do hold
similarly for the product of two Jacobi Triples, the dimension $d$ being replaced by
the sum $d + d'$ of the respective dimensions.

We quote finally the contraction property in the Kantorovich metric $W_2$
\begin {equation} \label {contraction}
W_2(P_t^*\nu, P_t^*\mu) \leq e^{-Kt} W_2(\nu, \mu),
\end {equation}
where $P_t^*$ is the heat semigroup acting (by duality) on measures,
as a consequence again of a curvature-dimension condition $CD(K,\infty)$ (cf.~\cite {BGL14}).

Before proceeding to the proof of the main results in the next sections, we mention here
the corresponding properties in the non-symmetric case with parameters $\alpha, \beta \geq \frac 12$.
Actually, all the above quantitative properties still hold for the non-symmetric Jacobi model with some
minor modifications. Particularly, the heat kernel bounds developed in \cite {NS13} imply that
\beq \label {eq.max}
p_t(x, x) \leq Ct^{- \max \{\alpha, \beta \} }, \quad x \in [-1, 1], \,\, 0 < t \leq 1.
\eeq
This conclusion also follows from a suitable Sobolev inequality as developed in \cite {B96}.
On the other hand, when $\alpha, \beta \geq \frac{1}{2}$, the curvature-dimension condition $CD(\frac{\alpha+\beta-1}{2} , 2(\alpha+\beta)-1)$ holds. Note that this condition does
not coincide with $CD(d-1,d)$ when $\alpha =\beta = \frac d2$ (it is weaker), but is good
enough to ensure all the curvature lower bounds necessary for the preceding semigroup properties.
Taking these properties for granted, the arguments developed below in the symmetric case extend
similarly to the non-symmetric Jacobi model.

\bigskip

To conclude this section, it is meaningful to make a comparison between the main results
emphasized in the introduction and what is known on the sphere, or the unit interval equipped with the
Euclidean metric. Recall the projection
operator
$$
T : x =(x_1, \ldots, x_{d+1}) \in \sph^d \mapsto x_1 \in [-1, 1]
$$
and denote by $\mu_{\sph^d}$ the uniform distribution on $\sph^d$, and by $X$
a random variable with distribution $\mu_{\sph^d}$.
Then it holds true that, with the obvious notation,
$$
 W_2^2( \mu^n, \mu ) \leq W_2^2( \mu_{\sph^d}^n, \mu_{\sph^d} ).
$$
To prove this assertion, notice that for any coupling $\pi$ on $\sph^d \times \sph^d$ with
marginals $\mu^n_{\sph^d} = \frac{1}{n} \sum_{i = 1}^{n} \delta_{X_i}$ and
$\mu_{\sph^d}$, the push-forward of $\pi$ by $T \otimes T$
is a coupling on $[-1, 1] \times [-1, 1]$ with marginals $\mu^n = \frac{1}{n} \sum_{i = 1}^{n} \delta_{T(X_i)}$ and $\mu$. Then,
using the spherical coordinate system, it is easy to see that for all $x, y \in \sph^d$,
$$
\rho( T(x), T(y) ) \leq \rho_{\sph^d} (x, y)
$$
where $\rho_{\sph^d}$ is the geodesic distance on $\sph^d$. Therefore by the definition of $W_2$,
the claim is proved since $T(X_i)$ has distribution $\mu$.

The known upper bounds for the sphere model (see e.g.~\cite{L17})
\begin{align*}
\E \left[ W_2^2( \mu_{\sph^d}^n, \mu_{\sph^d} ) \right] \lesssim \left\{ \begin{array}{cc}
\frac{1}{n},  & d = 1, \\
\frac{\log n}{n},  & d = 2, \\
\frac{1}{n^{2/d}}, & d \geq 3,
\end{array}
\right.
\end{align*}
thus transfer to the Jacobi model, at least when $d$ is an integer. However, unless $d=1$,
the rate given by Theorem~\ref {t11} is smaller than the rates on the sphere $\sph^d$.

\begin {remark}
It may be observed that the PDE proof of the lower bound in the Ajtai-Koml\'os-Tusn\'ady
provided in Section~5.2 of \cite {AST19} may be adapted, together with the suitable
version of Lusin's approximation of Sobolev functions, to show that the
rate $\frac{1}{n^{2/d}}$ is optimal on $\sph^d$ with $d \geq 3$ (and similarly on compact
manifolds satisfying the volume doubling property).
\end {remark}

In another direction, note that for $x, y \in [-1, 1]$, $\rho(x, y) \geq |x - y|$. Hence by the definition of $W_2$,
$ W_2^2( \mu^n, \mu ) \geq \widetilde{W}_2^2( \mu^n, \mu ) $, where $\widetilde{W}_2$ is
the Kantorovich distance defined via the standard metric $| \cdot |$ on the line.
Combining the relevant result obtained by ordered statistics \eqref {jacobieuclidean},
the order of growth of $\E [ W_2^2( \mu^n, \mu )]$ is at least
$\frac{1}{n}$, which is compatible with Theorem~\ref{t11}.

\section{Density fluctuation bounds } \label{S3}
\setcounter{equation}{0}

With this section, we start addressing the proofs of the main results, following the approach
of \cite {AST19}. The following objects and notation are taken from \cite{AST19}.
For $t \geq 0$, define
$$
r^n = \sqrt{n}(\mu^n - \mu), \quad \mu^{n, t} = P_t^* \mu^{n},  \quad
r^{n, t}\mu = P_t^* r^n = \sqrt{n}(\mu^{n, t} - \mu).
$$
In particular therefore
$$
r^{n, t}(y) = \int_E ( p_t(\cdot , y) - 1 ) dr^n
$$
with $ y \in E = [-1,1]$ or $[-1,1]^2$.

The heuristics here is that the empirical measure $\mu_n$ converges to $\mu$ as $n \to \infty$ in both
the weak and $W_2$ topologies (see e.g.~\cite{BL16}). On the other hand,
by definition of the heat kernel regularization,
for fixed $n$, $\mu^{n, t}$ is close to $\mu^{n}$ as $ t \to 0$.
The choice of the normalization of $r^n$ derives from the central limit theorem, and it is therefore
reasonable that if $t = t(n) \to 0$ is chosen properly, $\frac{r^{n, t}}{\sqrt{n}}$ will tend
to 0 as $n \to \infty$.

The following probabilistic statements quantify this heuristics, and
describe the relationship between the rate of $\frac{r^{n, t}}{\sqrt{n}}$ converging to 0 and $n$, $t$
in the probability sense. They are immediate versions of the corresponding
Propositions 3.9 and 3.10 of \cite{AST19}, following the analogous semigroup properties
described in the previous section. Uniformity in $y$ is achieved as in \cite {AST19} by
a covering argument. Namely, for
$\delta > 0$, if $N_E(\delta)$ denotes the smallest cardinality of a $\delta$-net of $E$,
then for $E=[-1, 1]$, $N_E(\delta) \lesssim \max \{1, \delta^{-1} \}$
while for $E=[-1, 1]^2$, $N_E(\delta) \lesssim \max \{1, \delta^{-2} \}$.

\begin{proposition} [Fluctuation bounds for the Jacobi model with $d \geq 1$] \label{p33}
There exist constants $C, C' > 0$ with the following property:
for all $\eta \in (0, 1)$, $\gamma \in (0, \frac{2}{d})$ and $\frac{1}{\eta n^\gamma} \leq t \leq C'$, we have
$$
\PP \bigg( \sup \limits_{y \in [-1, 1]} \frac{ |r^{n, t}(y)|}{\sqrt{n}} > \eta \bigg)
		 \leq C\exp \big( - \theta n^{1-\frac{d \gamma}{2}} \big)
$$
with $\theta = \theta (\eta) > 0$.
\end{proposition}

In parallel, the following proposition is valid on the product Jacobi model with $d = d' = 1$.
It is in particular in this result that the constraint $d+d' \leq 2$ appears, restricting the scope of
Theorem~\ref {t13}. Indeed, the product Jacobi model is of curvature-dimension
$CD(\min (d,d')-1, d+d')$ yielding the ultracontractivity property \textrm {(UC)} with
decay $t^{-(d+d')/2}$ for the small values of $t >0$. According to
Proposition 3.10 in \cite{AST19}, when $t$ if of the order of $\frac {\log n}{n}$,
the probabilistic tail estimate requires $d+d' \leq 2$. On the other hand,
$d,d' \geq 1$ is needed for the suitable heat kernel bounds.

\begin{proposition} [Fluctuation bounds for the product Jacobi model] \label{p34}
On the product Jacobi model $E = [-1, 1]^2$ with $d = d' = 1$, there exist constants $C, C'> 0$ with the following property: for all $\eta \in (0, 1)$ and $C' \geq t \geq \frac{\gamma \log n}{n}$
where $\gamma = \gamma(\eta) > 0$ is sufficiently large, we have
$$
\PP \left( \sup \limits_{y \in E} \frac{ |r^{n, t}(y)|}{\sqrt{n}} > \eta \right) \leq \frac{C}{n^2} .
$$
\end{proposition}

\section{Energy estimates } \label{S4}
\setcounter{equation}{0}

This section is devoted to the energy or trace estimates. It is here that the one-dimensional feature
of the Jacobi model is reflected, in contrast with the case of a compact manifold such as the sphere.

Recall that the eigenvalues of the $d$-dimensional Jacobi operator
$-L$ on $[-1,1]$ are given by the sequence $\lambda_k = k(k+d-1)$, $k \in \N$.
By the spectral representation of the heat kernel $p_s(x,y)$, $s >0$, $x, y \in [-1,1]$,
\begin{align*}
p_s(x, y) = \sum_{k = 0}^{\infty} e^{-sk(k+d-1)} J_k(x)J_k(y) ,
\end{align*}
and the trace formula may be expressed in the form
\begin{align} \label{tf}
\int_{-1}^{1} p_s(x, x) d\mu(x) = \sum_{k=0}^{\infty} e^{-sk(k+d-1)}
\end{align}
which will be essential for the validity of the next statement.

\begin{lemma}[Asymptotics for the trace of the Jacobi model with $d \geq 1$] \label{l41}
\begin{align*}
\int_{-1}^{1} p_s (x, x) d\mu(x)
   = \frac{1}{\sqrt{s}} \left( \frac{\sqrt{\pi}}{2} + o\left( 1 \right) \right) \quad \textrm{as }  s \to 0.
\end{align*}
\end{lemma}

\begin{proof}
Let $I(s) = \int_{0}^{\infty} e^{-su(u+d-1) } du$, $s >0$. Since
$$
e^{-s(k+1)(k+d)}  \leq e^{-su(u+d-1)} \leq e^{-sk(k+d-1)}
$$
for $k \leq u \leq k+1$, it follows that
\begin{align} \label{id}
I(s) \leq \sum_{k=0}^{\infty} e^{-sk(k+d-1)} \leq I(s) + 1.
\end{align}
Now
\begin{align*}
\sqrt{s}I(s) = \sqrt{s} \int_{0}^{\infty} e^{-su(u+d-1) } du = e^{\frac{(d-1)^2s}{4}} \sqrt{s} \int_{0}^{\infty} e^{-s(u+\frac{d-1}{2})^2 } du = e^{\frac{(d-1)^2s}{4}} \int_{\frac{(d-1)\sqrt{s}}{2}}^{\infty} e^{-t^2} dt
\end{align*}
so that
\begin{align*}
\lim_{s \to 0} \sqrt{s}I(s) = \int_{0}^{\infty} e^{-t^2} dt = \frac{\sqrt{\pi}}{2}.
\end{align*}
Combining with \eqref{id} and \eqref {tf}, the estimate is proved.
\end {proof}

On a product, the previous trace estimate immediately yields the following conclusion.

\begin{corollary} [Asymptotics for the trace of the product Jacobi model] \label{c42}
For the product Jacobi model $ E = [-1, 1]^k$ with ${d_j \geq 1}$, $j = 1, \ldots, k$,
$$
\int_E p_s(x, x) d\mu(x)
= \frac{1}{s^{\frac{k}{2}}} \left( \frac{\pi^{\frac{k}{2}}}{2^k} + o\left( 1 \right) \right) \quad \textrm{as }  s \to 0.
$$
\end{corollary}

\begin {remark}
It should be noted that the order of the trace as $s \to 0$ is independent of $d$ (and $d_j$).
\end {remark}

Let $f^{n, t}$ be the solutions of the PDE $Lf^{n, t} = r^{n, t} $ in $E=[-1, 1]$, whose means are zero.
The following results are consequences of the preceding trace asymptotics.
Again, the arguments follow the investigation \cite {AST19}.

\begin{lemma} [Energy estimates for the Jacobi model with $d \geq 1$] \label{l44}
If $t = t(n) \to 0$ as $n \to \infty$, then
\begin{align} \label{id1}
\lim_{n \to \infty} \E \left[ \int_{-1}^{1} \Gamma (f^{n, t}) d\mu \right]
= \int_{0}^{\infty} \left( \int_{-1}^{1} p_s(x, x) d\mu(x) -1 \right)ds
= \sum_{k=1}^{\infty} \frac{1}{k(k+d-1)}.
\end{align}
Moveover, if $1 \leq d < \frac{4}{3}$,
\begin{align} \label{id2}
\limsup_{n \to \infty} \E \left[ \left( \int_{-1}^{1} \Gamma (f^{n, t}) d\mu \right)^2 \right] < \infty .
\end{align}
When $ d \geq \frac{4}{3}$, if $t = t(n) \to 0$ as $n \to \infty$ and $t \geq \frac{c}{n^\gamma}$
for some $c > 0$, where $\gamma \in (0, \frac{2}{3d - 4})$, then \eqref {id2} also holds.
\end{lemma}

\begin {proof}
The proof of \eqref{id1} is an easy consequence of the trace formula,
and can be found in \cite{AST19}. We only prove \eqref{id2}. Using Lemma~3.16 of \cite{AST19},
\begin{align*}
\E \left[ \left( \int_{-1}^{1} \Gamma (f^{n, t}) d\mu \right)^2 \right]
&= \E \left[ \left( 2 \int_{t}^{\infty} \int_{-1}^{1} (P_sr^n)^2 d\mu ds \right)^2 \right]\\
& \leq 3\left[ \int_{2t}^{\infty} \left( \int_{-1}^1 p_s(x, x) d\mu(x) -1 \right) ds \right]^2 + \frac{1}{n} \left( 2 \int_{t}^{\infty} \int_{-1}^1 [\![p_s(\cdot, y)]\!]^2_4 d\mu(y) ds \right)^2 \\
& = 3\left( \sum_{k=1}^{\infty} \frac{1}{k(k+d-1)} \right)^2 + \frac{1}{n} \left( 2 \int_{t}^{\infty} \int_{-1}^{1} [\![p_s(\cdot, y)]\!]^2_4 d\mu(y) ds \right)^2,
\end{align*}
where $[\![p_s(\cdot, y)]\!]_4 = (\int_{-1}^{1} (p_s(x, y) - 1)^4 d\mu(x))^{\frac{1}{4}}$.

It is sufficient to estimate the limsup of the second term. Using the properties
(UC) and (SG), we know that for any $y \in [-1, 1]$,
\begin{align*}
[\![p_s(\cdot, y)]\!]^4_4  \lesssim \left\{ \begin{array}{cc}
s^{-\frac{3d}{2}},  &s \in (0, 1), \\
e^{-2Cs}, & s \in (1, \infty).
\end{array}
\right.
\end{align*}
Thus for $t \in (0, 1)$,
\begin{align*}
\int_{t}^{\infty} [\![p_s(\cdot, y)]\!]^2_4 ds \lesssim \left\{ \begin{array}{cc}
1,  &1 \leq d < \frac{4}{3}, \\
\log(\frac{1}{t}),  &d = \frac{4}{3}, \\
t^{-\frac{3d}{4}+1}, & d > \frac{4}{3}.
\end{array}
\right.
\end{align*}
Then it is easy to conclude. For example, if $ d > \frac{4}{3}$ and
$t \geq \frac{c}{n^\gamma}$ for some $c > 0$,
\begin{align*}
\frac{1}{n} \left( 2 \int_{t}^{\infty} \int_{-1}^{1} [\![p_s(\cdot, y)]\!]^2_4 d\mu(y) ds \right)^2 \leq \frac{C}{n^{1 - (\frac{3d}{2}-2)\gamma}},
\end{align*}
which is bounded if  $\gamma \in (0, \frac{2}{3d - 4})$. The proof of the lemma is complete.
\end {proof}

In the same way, the following energy estimates hold true for the product Jacobi model with $d = d' =1$.

\begin{lemma} [Energy estimates for the product Jacobi model] \label{l45}
On the product Jacobi model $E = [-1, 1]^2$ with $d = d' = 1$, if $t = t(n) \to 0$ as $n \to \infty$, then
\begin{align} \label{id4}
\lim_{n \to \infty} \frac{1}{|\log t| } \E \left[ \int_{E} \Gamma (f^{n, t}) d\mu \right] = \frac{\pi}{4}.
\end{align}
Furthermore, assuming that $t = t(n) > \frac{C}{n}$ for some $C>0$, we have
\begin{align} \label{id5}
\limsup_{n \to \infty} \frac{1}{(\log t)^2} \E \left[ \left( \int_{E} \Gamma (f^{n, t}) d\mu \right)^2 \right] < \infty .
\end{align}
\end{lemma}

\begin {proof}
Recall that
\begin{align*}
\E \left[ \int_{E} \Gamma (f^{n, t}) d\mu \right] = \int_{2t}^{\infty} \left( \int_{E} p_s(x, x) d\mu(x) -1 \right) ds.
\end{align*}
Combining Corollary \ref{c42} and (SG),
\begin{align*}
\int_{E} p_s(x, x) d\mu(x) - 1  \lesssim \left\{ \begin{array}{cc}
s^{-1},  &s \in (0, 1), \\
e^{-Cs}, & s \in (1, \infty).
\end{array}
\right.
\end{align*}
Thus, we obtain
\begin{align*}
\E \left[ \int_{E} \Gamma (f^{n, t}) d\mu \right] \lesssim | \log t | + 1,
\end{align*}
from which
$$
\limsup_{n \to \infty} \frac{1}{| \log t | } \E \left[ \int_{E} \Gamma (f^{n, t}) d\mu \right] < \infty .
$$
But what is actually really shown here is that
\begin{align*}
\frac{1}{|\log t| } \E \left[ \int_{E} \Gamma (f^{n, t}) d\mu \right] - \frac{1}{|\log t| } \int_{2t}^{1} \left( \int_{E} p_s(x, x) d\mu(x) -1 \right) ds
\end{align*}
is infinitesimal as $ n \to \infty$. Combining this with Corollary \ref{c42}, we obtain \eqref{id4}.

Now we turn to \eqref{id5}. Recall, as in the previous lemma, that
\begin{align*}
\E \left[ \left( \int_{E} \Gamma (f^{n, t}) d\mu \right)^2 \right]
&= \E \left[ \left( 2 \int_{t}^{\infty} \int_{E} (P_sr^n)^2 d\mu ds \right)^2 \right]\\
& \leq 3\left[ \int_{2t}^{\infty} \left( \int_{E} p_s(x, x) d\mu(x) -1 \right) ds \right]^2 + \frac{1}{n} \left( 2 \int_{t}^{\infty} \int_{E} [\![p_s(\cdot, y)]\!]^2_4 d\mu(y) ds \right)^2 .
\end{align*}
By the preceding,
\begin{align*}
\int_{2t}^{\infty} \left( \int_{E} p_s(x, x) d\mu(x) -1 \right) ds \lesssim |\log t|.
\end{align*}
Thus it is sufficient to estimate the limsup of the second term. Using the properties of $p_t$,
and again that $d=d'=1$, we know that for any $y \in E$,
\begin{align*}
[\![p_s(\cdot, y)]\!]^4_4  \lesssim \left\{ \begin{array}{cc}
s^{-3},  &s \in (0, 1), \\
e^{-Cs}, & s \in (1, \infty).
\end{array}
\right.
\end{align*}
Then it is easy to verify that if $t = t(n) > \frac{C}{n}$ for some $C>0$,
\begin{align*}
\frac{1}{n} \left( 2 \int_{t}^{\infty} \int_{-1}^{1} [\![p_s(\cdot, y)]\!]^2_4 d\mu(y) ds \right)^2 \lesssim \frac{1}{nt},
\end{align*}
which is bounded. The proof is therefore complete.
\end {proof}

\section{Regularization estimates} \label{S5}
\setcounter{equation}{0}

In this section, we evaluate the cost of the regularization by the heat kernel. We make use
here of some recent advances from \cite {AG18} which provide stronger regularization error bounds
when dealing with empirical measures.

To start with, we state a general result from \cite {AST19} (see also \cite {AG18,L17})
on the connection between the solution of the Poisson equation
$Lf = u$ and the Kantorovich distance. The original proof is
developed in the setting of Riemannian manifolds, but works similarly here.
The following theorem is stated for $E = [-1,1]$ and $E = [-1,1]^2$.
The use of the integration by parts formula in the proof requires
to add some Neumann boundary conditions to the solution of the Poisson equation
$Lf = u$. But for example, for 1-dimensional Jacobi models, $\Gamma(f)(x) = (1- x^2)f'^2 = 0$ holds naturally when $x = \pm 1$ for all smooth $f$.

\begin{theorem} \label{t43}
Given smooth, positive density functions $u_0$, $u_1$ in
$E$, let $f$ be the unique solution of $Lf = u_0 - u_1$ with mean zero. Then
\begin{align*}
W_2^2( u_0\mu, u_1\mu) \leq 4 \int_E \frac{\Gamma(f)}{u_0} d\mu
\end{align*}
and
\begin{align*}
W_2^2( u_0\mu, u_1\mu) \leq \int_E \Gamma(f) \frac{\log(u_1) - \log(u_0)}{u_1 - u_0} d\mu.
\end{align*}
\end{theorem}

In the following, we start the study of the regularization procedure
on $E = [-1,1]$. First, since $\mu= \int_{-1}^{1} \delta_x d\mu(x)$ and
$$
W_2^2( P_t^*\delta_x, \delta_x ) \leq \int_{-1}^{1} \rho^2(x, y) p_t(x, y) d\mu(y) ,
$$
from the joint convexity of $W_2^2$ and (DR), for every probability $\mu $ (on $E = [-1,1]$),
\begin{align*}
\E \left[W_2^2( P_t^* \mu, \mu )\right]
\leq \int_{-1}^{1} \int_{-1}^{1} \rho^2(x, y) p_t(x, y) d\mu(x) d\mu (y) \leq Ct.
\end{align*}
Hence
\begin {equation} \label {reg}
 \E \left[W_2^2( \mu^n, \mu^{n, t} )\right] \leq Ct.
\end {equation}
This estimate is good enough for the case $1 \leq d < 2$. But for the case $d \geq 2$, we need the following more refined estimate. The idea is based on Theorem~5.2 of the recent \cite{AG18}.

\begin{theorem} \label{t51}
 Given $n \geq 1$, let the event
\begin{align*}
A_{\frac{1}{2}, n}
  = \bigg \{ \sup \limits_{y \in [-1, 1]} \frac{ |r^{n, t}(y)|}{\sqrt{n}} \leq \frac{1}{2} \bigg \}.
\end{align*}
Then for the Jacobi model with $d \geq 1$, for $t \in (0, 1)$, we have
\begin{align*}
\E \left[ W_2^2( \mu^n, \mu^{n, t} ) \right] \lesssim \PP ( A_{\frac{1}{2}, n}^c ) + \frac{\sqrt{t}}{n}.
\end{align*}
\end{theorem}

\begin {proof}
Fix a time $t_0 \in (0, t)$. Applying \eqref {reg} at time $t_0$, it follows from the triangle inequality that
\begin{align*}
\E \left[ W_2^2( \mu^n, \mu^{n, t} ) \right]  &\leq 2\E \left[ W_2^2( \mu^n, \mu^{n, t_0} )\right] + 2\E \left[ W_2^2( \mu^{n, t_0}, \mu^{n, t}  ) \right]\\
& \lesssim t_0 + \E \left[ W_2^2( \mu^{n, t_0}, \mu^{n, t}  ) \right].
\end{align*}
In order to estimate $\E [W_2^2( \mu^{n, t_0}, \mu^{n, t}  )]$, let $f$ be the solution of $Lf = u^{n, t} - u^{n, t_0}$, where $u^{n, t_0}$ and $u^{n, t}$ are the density functions of $\mu^{n, t_0}$ and $\mu^{n, t}$ with respect to $\mu$, respectively. By definition,
\begin{align*}
|u^{n, t}-1| = \frac{ |r^{n, t}(y)|}{\sqrt{n}}.
\end{align*}
Therefore in the event $A_{\frac{1}{2}, n}$, using Theorem~\ref{t43}, we obtain
\begin{align*}
W_2^2( \mu^{n, t_0}, \mu^{n, t}  ) \leq 4\int_{-1}^{1} \frac{\Gamma(f)}{u^{n, t}} d\mu \leq 8\int_{-1}^{1} \Gamma(f) d\mu.
\end{align*}
Hence, due to the independence of the variables $X_i$,
\begin{align*}
\E \left[ W_2^2( \mu^{n, t_0}, \mu^{n, t}  ) \chi_{A_{\frac{1}{2}, n}} \right]
& \lesssim \E \left[ \int_{-1}^{1} \Gamma(f) d\mu \right] \\
&= \E \left[ \int_{-1}^{1} (-Lf)  f d\mu \right]\\
& = \E \left[ \int_{-1}^{1} (u^{n, t_0} - u^{n, t}) \left( \int_{t_0}^{t} u^{n, s} ds \right) d\mu \right]\\
& = \frac{1}{n} \int_{t_0}^{t} \E \left[ \int_{-1}^{1} (p_{t_0}(X, y) - p_t(X, y))p_s(X, y) d\mu(y) \right] ds\\
& = \frac{1}{n} \int_{t_0}^{t} \int_{-1}^{1} (p_{t_0+s}(x, x) - p_{t+s} (x, x) ) d\mu(x) ds\\
&\lesssim \frac{1}{n} \int_{2t_0}^{t+t_0} \int_{-1}^{1} (p_s(x, x) - 1) d\mu(x) ds.
\end{align*}
Lemma~\ref{l41} tells us that
$$
\E \left[ W_2^2( \mu^{n, t_0}, \mu^{n, t}  ) \chi_{A_{\frac{1}{2}, n}} \right]
\lesssim \frac{1}{n} \int_{2t_0}^{t+t_0}  \frac{1}{\sqrt{s}} ds
 \lesssim \frac{1}{n} \frac{t-t_0}{\sqrt{t}}.
$$
On the other hand $W_2^2( \mu^n, \mu ) \leq \pi^2$. Therefore
\begin{align*}
\E \left[ W_2^2( \mu^n, \mu^{n, t} ) \right]
&\lesssim t_0 + \E \Big[ W_2^2( \mu^{n, t_0}, \mu^{n, t}  ) \chi_{A_{\frac{1}{2}, n}^c} \Big]
+ \E \Big[ W_2^2( \mu^{n, t_0}, \mu^{n, t}  ) \chi_{A_{\frac{1}{2}, n}} \Big]\\
&\lesssim t_0 + \PP ( A_{\frac{1}{2}, n}^c ) + \frac{1}{n} \frac{t-t_0}{\sqrt{t}}.
\end{align*}
Letting $t_0 \to 0$, the theorem is proved.
\end {proof}

If we choose $t=t(n) \to 0$ as $n \to \infty$ properly, then Proposition~\ref{p33} implies that the term $\PP ( A_{\frac{1}{2}, n}^c )$ is of much lower order than $\frac{\sqrt{t}}{n}$ . Therefore we obtain

\begin{corollary} [Refined regularization estimate for the Jacobi model with $d \geq 1$] \label{c52}
Assume that $ d \geq 1$. If $t = t(n) \to 0$ as $n \to \infty$ and $t \geq \frac{c}{n^\gamma}$ for some $c > 0$, where $\gamma \in (0, \frac{2}{d})$,
\begin{align*}
\lim_{n \to \infty} n \E \left[ W_2^2( \mu^n, \mu^{n, t} ) \right] = 0.
\end{align*}
\end{corollary}

With the same arguments, and Proposition~\ref{p34},
we deduce a corresponding result for the product Jacobi model with $d = d' =1$.

\begin{theorem} [Refined regularization estimate for the product Jacobi model] \label{t53}
On the product Jacobi model $E = [-1, 1]^2$ with $d = d' = 1$, there exists a constant
$C > 0$, such that for $t = t(n) = \frac{\gamma(n)}{n} \to 0$ with $\gamma(n) \geq C \log n$
\begin{align*}
\E \left[W_2^2( \mu^n, \mu^{n, t} )\right] \lesssim \frac{\log \gamma(n) }{n}.
\end{align*}

\end{theorem}

\section{Proofs of Theorems~\ref{t11} and \ref {t13}}\label{S6}
\setcounter{equation}{0}

We first address the proof of Theorem~\ref{t11}, which is split in an upper bound
and a lower bound. The upper bound, developed in the next Proposition~\ref {p61},
simplifies the corresponding step in \cite {AST19}
via the improved regularization Theorem~\ref {t51} and Corollary~\ref {c52}
(and in particular avoids the use of Riesz transform bounds).

\begin{proposition} \label{p61}
 Assume that $d \geq 1$, then
\begin{align*}
\limsup_{n \to \infty} n\E \left[ W_2^2( \mu^n, \mu ) \right] \leq \sum_{k=1}^{\infty} \frac{1}{k(k+d-1)}.
\end{align*}
\end{proposition}

\begin {proof}
Fix $\gamma \in (0, \frac{2}{d})$ and let $t = t(n) =  n^{-\gamma}$. For $\eta \in (0, 1)$ consider the event
\begin{align*}
A_\eta = A_{\eta, n} = \bigg \{ \sup \limits_{y \in [-1, 1]} \frac{ |r^{n, t}(y)|}{\sqrt{n}} \leq \eta \bigg\}.
\end{align*}
Since $W_2^2( \mu^n, \mu ) \leq \pi^2$, for $n$ large enough using Proposition~\ref{p33},  we have
\begin{align*}
n\E \left[ W_2^2( \mu^n, \mu ) \right] &= n\E \left[ W_2^2( \mu^n, \mu )\chi_{A_\eta} \right] + n\E \left[ W_2^2( \mu^n, \mu )\chi_{A_\eta^c} \right]\\
& \leq n\E \left[ W_2^2( \mu^n, \mu )\chi_{A_\eta} \right] + Cn\exp \left( - \theta n^{1-\frac{d \gamma}{2}} \right)
\end{align*}
with $\theta = \theta (\eta) > 0$.
Therefore it is sufficient to estimate $n\E [ W_2^2( \mu^n, \mu ) \chi_{A_\eta} ]$.

Notice that for $\varepsilon > 0$,
$$
W_2^2( \mu^n, \mu ) \leq (1 + \varepsilon) W_2^2( \mu^{n, t}, \mu ) + (1 + \varepsilon^{-1}) W_2^2( \mu^n, \mu^{n, t} ).
$$
Corollary~\ref{c52} tells us that $\lim_{n \to \infty} n \E \left[ W_2^2( \mu^n, \mu^{n, t} ) \right] = 0$.
So we only need to estimate
\begin{align*}
\limsup_{n \to \infty} n\E \left[ W_2^2( \mu^{n, t}, \mu )\chi_{A_\eta} \right].
\end{align*}

On the event $A_\eta$, using Theorem~\ref{t43}, we obtain
\begin{align*}
W_2^2( \mu^{n, t}, \mu )
 \leq \frac{1}{n} \int_{-1}^{1} \Gamma(f^{n, t})\frac{\log(u^{n,t})}{u^{n, t} - 1}  d\mu
\leq \frac{1}{n\sqrt{1-\eta}} \int_{-1}^{1} \Gamma(f^{n, t}) d\mu,
\end{align*}
where $u^{n,t} = \frac{ r^{n, t}}{\sqrt{n}} $ is the density function of $\mu^{n, t}$ with respect to $\mu$
(since $\frac {\log u}{|u-1|} \leq \frac {1}{\sqrt {1 - \eta}}$ whenever $u >0$, $|u-1| \leq \eta$,
$0 < \eta < 1$).
Hence Lemma~\ref{l44} implies that
\begin{align*}
\limsup_{n \to \infty} n\E \left[ W_2^2( \mu^{n, t}, \mu )\chi_{A_\eta} \right] &\leq \lim_{n \to \infty} \frac{1}{\sqrt{1 - \eta}} \, \E \left[ \int_{-1}^{1} \Gamma (f^{n, t}) d\mu \right]\\
& = \frac{1}{\sqrt{1- \eta}} \sum_{k=1}^{\infty} \frac{1}{k(k+d-1)}.
\end{align*}
At last we get that
\begin{align*}
\limsup_{n \to \infty} n\E \left[ W_2^2( \mu^n, \mu ) \right] \leq \frac{1+\varepsilon}{\sqrt{1- \eta}} \sum_{k=1}^{\infty} \frac{1}{k(k+d-1)}.
\end{align*}
Letting $\varepsilon, \eta \to 0$, the proposition is established.
\end {proof}

Next, we prove the lower bound part of Theorem~\ref{t11} in the form of the following statement.

\begin{proposition}Assume that $d \geq 1$, then
\begin{align*}
\liminf_{n \to \infty} n\E \left[ W_2^2( \mu^n, \mu ) \right] \geq \sum_{k = 1}^{\infty} \frac{1}{k(k+d-1)}.
\end{align*}
\end{proposition}

\begin {proof}
First, by the heat flow contraction in Wasserstein space \eqref {contraction}, for every $t \geq 0$,
\begin{align*}
W_2^2( \mu^n, \mu ) \geq e^{2(d-1)t} \, W_2^2( \mu^{n, t} , \mu ) \geq W_2^2( \mu^{n, t} , \mu )
\end{align*}
so that we need only concentrate on $W_2^2( \mu^{n, t} , \mu )$.
Next, by the Kantorovich duality (cf.~\cite {V09}),
\begin{align*}
\frac{1}{2} W_2^2( \mu^{n, t} , \mu )
&\geq \sup \left\{ \int_{-1}^{1} f d\mu^{n, t} + \int_{-1}^{1} g d\mu
\; ; \; f(x)+g(y)\leq \frac{\rho(x,y)^2}{2} \right\}\\
& = \sup \left\{ \int_{-1}^{1} f \frac{r^{n, t}}{\sqrt{n}} + \int_{-1}^{1} (f + g) d\mu
\; ;\; f(x)+g(y)\leq \frac{\rho(x,y)^2}{2} \right\}.
\end{align*}

Fix $\eta \in (0, 1)$ and let $t = t(n) = \eta^{-1} n^{-\gamma}$, where
\begin{align*}
\gamma \in \left\{ \begin{array}{cc}
\left(0, \frac{2}{d} \right),  &1 \leq d \leq \frac{4}{3}, \\
\left(0, \min \{ \frac{2}{d}, \frac{2}{3d-4} \} \right) , & d > \frac{4}{3}.
\end{array}
\right.
\end{align*}
As in the proof of Proposition~\ref{p61}, consider the event $A_\eta$.
Define $f = -\frac{f^{n, t}}{\sqrt{n}}$, where $f^{n, t}$ is defined in Section \ref{S3}.
Therefore, on the event $A_\eta$, $\| Lf \|_{\infty} \leq \eta$.
Moreover, using \eqref{i1},
$$
\| Lf \|_{\infty} + \frac{e^{2(d - 1) t} - 1}{2} \, \| \Gamma(f) \|_{\infty} \leq \omega(\eta)
$$
with $\omega(\eta) \to 0$ as $\eta \to 0$. To this $f$, we associate
$$
g(x) = Q_1(-f)(x) = \inf_{y \in [-1, 1]} [-f(y) + \frac{1}{2}\rho^2(x, y)] ,
$$
so that on $A_\eta$, by Corollary~3.3 of \cite{AST19},
\begin{align*}
\frac{1}{2} W_2^2( \mu^{n, t} , \mu ) \geq \left( 1 - \frac{e^{\omega(\eta)}}{2} \right)\frac{1}{n} \int_{-1}^{1} \Gamma(f^{n, t}) d\mu.
\end{align*}
Next
\begin{align*}
\frac{1}{2-e^{\omega(\eta)}} \liminf_{n \to \infty} & n\E \left[ W_2^2( \mu^{n, t} , \mu )
				 \chi_{A_\eta} \right] \\
&\geq \liminf_{n \to \infty} \E \left[ \chi_{A_\eta} \int_{-1}^{1} \Gamma(f^{n, t}) d\mu \right]\\
&\geq \lim_{n \to \infty} \E \left[ \int_{-1}^{1} \Gamma (f^{n, t}) d\mu \right] - \limsup_{n \to \infty} \E \left[ \chi_{A_\eta^c} \int_{-1}^{1} \Gamma(f^{n, t}) d\mu \right] . \\
\end{align*}
By Lemma~\ref{l44}, for our choice of $t = t(n)$,
\begin{align*}
\limsup_{n \to \infty} \E \left[ \left( \int_{-1}^{1} \Gamma (f^{n, t}) d\mu \right)^2 \right] < \infty.
\end{align*}
Hence
\begin{align*}
\limsup_{n \to \infty} \E \left[ \chi_{A_\eta^c} \int_{-1}^{1} \Gamma(f^{n, t}) d\mu \right] \leq \limsup_{n \to \infty} \PP(A_{\eta,n}^c)^{1/2}\left( \E \left[ \left( \int_{-1}^{1} \Gamma (f^{n, t}) d\mu \right)^2 \right]\right)^{1/2} = 0.
\end{align*}
Therefore
$$
\frac{1}{2-e^{\omega(\eta)}} \liminf_{n \to \infty}  n\E \left[ W_2^2( \mu^{n, t} , \mu )
				 \chi_{A_\eta} \right]
		\geq \lim_{n \to \infty} \E \left[ \int_{-1}^{1} \Gamma (f^{n, t}) d\mu \right]
		=  \sum_{k=1}^{\infty} \frac{1}{k(k+d-1)}
$$
from which the proposition follows by letting $\eta \to 0$.
\end {proof}

The modifications to address in the same way Theorem~\ref {t13} are minor. For Theorem~\ref {t13}, Proposition~\ref {p33} is replaced by Proposition~\ref {p34}, and
Lemma~\ref {l44} by Lemma~\ref {l45}. We only need to take $t(n) = \frac{\gamma \log n}{n}$, where $\gamma > 0$ is a sufficiently large constant. Then, the regularization
procedure is here taken into account by Theorem~\ref{t53} which indicates that
\begin{align*}
\frac{n}{\log n} \E \left[W_2^2( \mu^n, \mu^{n, t} )\right] \lesssim \frac{\log \log n}{\log n}
\end{align*}
so that in particular
$\lim_{n \to \infty} \frac{n}{\log n} \E \left[ W_2^2( \mu^n, \mu^{n, t} ) \right] = 0$.

By the semigroup properties discussed in Section~\ref{S2},
the corresponding results of Proposition~\ref {p33}, Lemma~\ref {l44} and Corollary~\ref {c52}
are valid for non-symmetric Jacobi models, with the value of $d$ therein replaced by
$2(\alpha +\beta) - 1$ due to the curvature-dimension
$CD(  \frac {\alpha + \beta -1}{2}, 2(\alpha +\beta) - 1)$ condition.
All the other steps are then basically similar to the proof of Theorem~\ref {t11} and are thus
not repeated.

In the same way, the proofs of the limits
\eqref {bipartite1} and \eqref {bipartite2} in the bipartite case are
immediately adapted from Propositions~4.9 and 4.10 of \cite{AST19}.

\section{Proof of Theorem~\ref{t12}} \label{S7}
\setcounter{equation}{0}

We cover at the same time Theorem~\ref{t12} and its extension \eqref {eq.kjacobi} to $k \geq 3$.
Let therefore ${E = [-1, 1]^k}$ and $\mu$ be the product Jacobi measure with $d_j \geq 1$,
$j = 1, \ldots, k$ ($k \geq 2$). The arguments also apply to non-symmetric Jacobi measures
(with the conditions $\alpha_j, \beta_j \geq \frac 12$).

The following lemma is a general heat kernel upper bound of
$\E \left[ W_2^2( \mu^n, \mu ) \right]$ (for its proof, see \cite{L17})

\begin{lemma} \label{l81}
In the prescribed setting and notation, for every $t > 0$,
\begin{align*}
\E \left[ W_2^2( \mu^n, \mu ) \right] \leq 2 \int_E \int_E \rho^2(x, y) p_t(x, y) d\mu(x) d\mu(y) + \frac{8}{n}\int_{2t}^{\infty} \int_E (p_s(x, x) - 1) d\mu(x) ds.
\end{align*}
\end{lemma}

Combining (DR) and Corollary~\ref{c42},
we know that on the product Jacobi model, for every $t \in (0, 1)$,
\begin{align*}
\E \left[ W_2^2( \mu^n, \mu ) \right] \lesssim t + \left\{ \begin{array}{cc}
 \frac{\log \left( \frac{1}{t} \right) }{n},  &k = 2, \\
 \frac{1}{nt^{(k / 2) - 1}}, & k \geq 3.
\end{array}
\right.
\end{align*}
After optimization in $t \in (0, 1)$, we obtain that
\begin{align*}
\E \left[ W_2^2( \mu^n, \mu ) \right] \lesssim \left\{ \begin{array}{cc}
\frac{\log n}{n},  & k = 2, \\
\frac{1}{n^{2/k}}, & k \geq 3,
\end{array}
\right.
\end{align*}
which amounts to the first half of Theorem~\ref {t12}.

In order to prove the lower bound part of Theorem~\ref{t12}, we use a simple comparison argument. Denote by $\lambda$ the uniform distribution on $[0, 1]^k$, let $U_1, \ldots, U_n$ be independent and distributed as $\lambda$ and write $\nu_n = \frac{1}{n} \sum_{i=1}^{n} \delta_{U_i}$. Define functions
$\Phi_j : [-1, 1] \to [0, 1]$, $j = 1, \ldots, k$, as follows:
\begin{align*}
\Phi_j (x) = C_{d_j} \int_{-1}^{x} (1 - u^2)^{\frac{d_j}{2} - 1} du,
\end{align*}
where $C_{d_j}$ are normalization constants. Write
$\Phi = \Phi_1 \otimes \cdots \otimes \Phi_k : E \to [0, 1]^k$. Therefore if $X_1, \ldots, X_n$ are  independent distributed as $\mu$, then $\Phi(X_1), \ldots, \Phi(X_n)$ are independent distributed as $\lambda$. The Kantorovich distance $W_1$ in this context is given by
\begin{align} \label{id7}
W_1( \mu^n, \mu ) = \sup_{\Gamma(f) \leq 1}
  \bigg[ \frac{1}{n} \sum_{i=1}^{n} f(X_i) - \int_{E} f d\mu \bigg].
\end{align}
Notice that for any smooth function $g$ defined on $[0, 1]^k$,
$g \circ \Phi$ is a function defined on $E$ and
\begin{align*}
\Gamma( g \circ \Phi ) & = \sum_{j = 1}^{k} (1-x_j^2) \Phi_j^{' 2} (x_j) (\partial_{x_j} g)^2  = \sum_{j = 1}^{k} C_{d_j}^2 (1-x_j^2)^{d_j - 1} (\partial_{x_j} g)^2.
\end{align*}
Therefore if $d_j \geq 1$, $j = 1, \ldots, k$ (or
$\alpha_j, \beta_j \geq \frac{1}{2}$), there exists a constant $C > 0$ such that
\begin{align} \label{id8}
\Gamma(g \circ \Phi) \leq C^2 \| \nabla g \|_{\infty}^2.
\end{align}
Combining with \eqref{id7}, we obtain that
\begin{align*}
W_1( \mu^n, \mu ) &\geq \sup_{\Gamma(g \circ \Phi) \leq 1} \bigg[ \frac{1}{n} \sum_{i=1}^{n} g \circ \Phi(X_i) - \int_{E} g \circ \Phi d\mu \bigg]\\
&\geq \sup_{ \| \nabla g \|_{\infty} \leq \frac{1}{C}} \bigg[ \frac{1}{n} \sum_{i=1}^{n} g \circ \Phi(X_i) - \int_{E} g \circ \Phi d\mu \bigg]\\
& = \frac{1}{C} \sup_{ \| \nabla g \|_{\infty} \leq 1} \bigg[ \frac{1}{n} \sum_{i=1}^{n} g  (\Phi(X_i)) - \int_{[0, 1]^k} g  dx dy \bigg]\\
& = \frac{1}{C} W_1( \nu^n, \lambda ).
\end{align*}
Therefore there exists a constant $C >0$,
\begin{align*}
\E \left[ W_2^2( \mu^n, \mu ) \right] \geq \E \left[ W_1( \mu^n, \mu ) \right]^2
	\geq \frac 1C \E \left[ W_1( \nu^n, \lambda ) \right]^2.
\end{align*}
From \eqref {AKT},
\begin{align*}
\E[W_1( \nu^n, \lambda )] \approx \left\{ \begin{array}{cc}
\sqrt{\frac{\log n}{n}},  & k = 2, \\
\frac{1}{n^{1/k}}, & k \geq 3.
\end{array}
\right.
\end{align*}
Theorem~\ref{t12} is proved.

\vskip 7mm

\textit {Acknowledgements}. I thank M.~Ledoux for useful indications and comments
during the preparation of this work, and the referee for several helpful suggestions for improvements.

\vskip 8mm

\font\tenrm =cmr10  {\tenrm

\parskip 0mm

\noindent School of Mathematical Sciences

\noindent Fudan University, Shanghai 200433, People's Republic of China

\noindent 15110840006@fudan.edu.cn

\noindent and

\noindent Institut de Math\'ematiques de Toulouse

\noindent Universit\'e de Toulouse -- Paul-Sabatier, F-31062 Toulouse, France

\noindent zhu@math.univ-toulouse.fr

}

\end{document}